\documentclass[10pt,amssymb]{article}

\usepackage{amsfonts,amssymb,amsmath,latexsym,verbatim,amscd}
\usepackage{amsthm,amsmath,amssymb,amscd}
\usepackage{amssymb}
\usepackage{graphics}
\usepackage{fancyhdr}

\def\à{\`a}
\def\é{\'e}
\def\è{\`e}
\def\ç{\'o}
\def\ù{\`u}
\def\ò{\`o}
\def\ì{\`\i\ }

\newcommand{\vol}{{\rm vol}}
\newcommand{\di}{{\rm d}}
\newcommand{\dvol}[1]{{\rm \,d}vol_{#1}}

\newcommand{\nor}[2]{\|{#1}\|_{#2}}
\newcommand{\al}{\alpha}
\newcommand{\be}{\beta}
\newcommand{\bigo}[1]{\mathcal{O} \left( #1 \right)}
\newcommand{\littleo}[1]{\mathcal{\sigma} \left( #1 \right)}

\newcommand{\zi}{z^{i}}

\newcommand{\ep}{\varepsilon}
\newcommand{\ex}{{\rm{e}}}
\newcommand{\te}{\theta}

\newcommand{\gep}{g_{\ep}}

\newcommand{\uep}{u_{\ep}}
\newcommand{\tip}{{\rm {th}}}
\newcommand{\cip}{{\rm {ch}}}

\newcommand{\norg}[2]{\nor{#1}{\mathcal{C}^0_{\gamma}(#2)}}
\newcommand{\norgd}[2]{\nor{#1}{\mathcal{C}^0_{\gamma + 2}(#2)}}
\newcommand{\tup}{\tilde{u}_P}
\newcommand{\tf}{\tilde{f}}

\newcommand{\tu}{\tilde{u}}
\newcommand{\nogd}[1]{\nor{#1}{\mathcal{C}^0_{\gamma + 2}}}
\newcommand{\nog}[1]{\nor{#1}{\mathcal{C}^0_{\gamma}}}

\newtheorem{teor}{Theorem}[section]

\newtheorem{lemma}[teor]{Lemma}
\newtheorem{prop}[teor]{Proposition}

\title{Generalized connected sum construction for scalar flat metrics}

\author{L. Mazzieri~\thanks{Scuola Normale Superiore di
Pisa and Laboratoire d'Analyse et de Math\'ematiques Appliqu\'ees,
Universit\'e Paris 12. E-mail : l.mazzieri@sns.it}}

\begin{document}

\maketitle

\section{Introduction and statement of the result}

In this paper we will show that the generalized connected sum
construction for constant scalar curvature metrics can be extended
to the zero scalar curvature case. In particular we want to
construct solutions to the Yamabe equation on the generalized
connected sum $M = M_1 \, \sharp_K \, M_2$ of two compact Riemannian
manifolds $(M_1,g_1)$ and $(M_2,g_2)$ with zero constant scalar
curvature along a common (isometrically embedded) submanifold
$(K,g_K)$ of codimension $\geq 3$.\\

We present here two kinds of construction. The first one is the
basic model and it works for every couple of scalar flat manifolds,
but it has a drawback. In fact following this method we are not
allowed to choose a scalar flat metric on the generalized connected
sum, although the error can be chosen as small as we want. The
second construction is an adjustment of the first one which enable
us to get a zero scalar curvature metric on the final manifold, but
it require the hypothesis that the starting Riemannian manifolds are
non Ricci flat.\\

In section 2-5 we present the first method. As in the nonzero scalar
curvature case, our strategy lies in writing down a family of
approximate solution metrics $(g_{\ep})_{\ep \in (0,1)}$ (where the
parameter $\ep$ represents the size of the tubular neighborhood we
excise from each manifold in order to perform the generalized
connected sum) and then in finding out a conformal factor $u_{\ep}$
such that for sufficiently small $\ep > 0$ the metrics
$\tilde{g}_{\ep} = u_\ep^{\frac{4}{m-2}} g_\ep$, $\ep \in (0,1)$,
are "small" constant scalar curvature metrics. As we claimed before,
notice that by this method it is impossible to ensure that the
scalar curvature $S = S_{\tilde{g}_{\ep}}$ of the metrics we obtain
is exactly zero. Anyway we will show that $S = \bigo{\ep^{n-2}}$.
Notice also that in order to achieve our goal we will need to scale
up or down the initial manifolds by means of suitable homotheties,
in other words we need to multiply the initial metrics $g_1$ and
$g_2$ by suitable positive constants; hence, what the submanifold
$(K,g_K)$ is actually required to do is to be homothetically
embedded in both the initial manifolds, it is to say isometrically
embedded modulo homotheties in $M_1$ and $M_2$. Let us
now describe this result more precisely.\\

Let $(M_1, g_1)$ and $(M_2,g_2)$ be two $m$-dimensional compact
Riemannian manifolds with zero constant scalar curvature, and
suppose that there exists a $k$-dimensional Riemannian manifold
$(K,g_K)$ which is isometrically embedded in each $(M_i,g_i)$, for
$i=1,2$, $m \geq 3$, $m-k \geq 3$. We also assume that the normal
bundles of $K$ in $(M_i,g_i)$ can be diffeomorphically identified.\\

Let $M^{R,Q} = M_1 \,\sharp_{K}\, M_2$ be the generalized connected
sum of $(M_1, R g_1)$ and $(M_2, Q g_2)$ along $K$ which is obtained
by removing an $\ep$-tubular neighborhood of $K$ from each $M_i$ and
identifying the two boundaries.\\

Our main result reads~:
\begin{teor}
Under the above assumptions, for $\ep \in (0, \ep_0)$ and suitable
constants $R,Q>0$, it is possible to endow $M^{R,Q}$ with a family
of constant scalar curvature metrics $(\tilde{g}_{\ep})$, whose
scalar curvature $S_{\tilde{g}_{\ep}}$ is a $\bigo{\ep^{n-2}}$. In
addition the metric $\tilde{g}_{\ep}$ is conformal to the metrics
$g_i$ away from a fixed (small) tubular neighborhood of $K$ in
$M_i$, $i=1,2$ for a conformal factor $u_\ep$ which can be chosen so
that
\[
\nor{u_{\ep} - 1}{L^{\infty}(M)} \leq C r_\ep
\]
where $r_\ep = \bigo{\ep^{1-\gamma}}$, $\gamma \in (0,1)$, for $n=3$
and $r_\ep = \bigo{\ep}$ for $n \geq 4$.
\end{teor}

Section 6 is devoted to the description of a special device, which
works in the non Ricci flat case. In this case we will be able to
achieve a scalar flat metric on the final manifold. The strategy
lies in making a slight modification of the approximate solution
metrics away from the polyneck. If the starting manifolds are non
Ricci flat, this construction provide us two correction terms which
will be employed in the nonlinear analysis in place of the non zero
constant scalar curvature and in place of the homotheties in order
to get a solution of the Yamabe equation with prescribed zero scalar
curvature.\\

The statement of the theorem is the following
\begin{teor}
Let $M$ be the generalized connected sum of two Riemannian scalar
flat non Ricci flat manifolds $(M_1,g_1)$ and $(M_2,g_2)$ of
dimension $m \geq 3$ along a common isometrically embedded
submanifold $(K, g_K)$ of codimension at least $3$. Under these
assumptions it is possible to endow $M$ with a family of scalar flat
metrics
\end{teor}

\section{Geometric construction}

The geometric construction we use here is essentially the same we
used in \cite{Mazz}, but in order to fix the notation it is useful
to transfer it, paying attention in making the appropriate
adjustments. Here we describe the construction in the case where $R
= Q = 1$, but it still holds in the general case. Of course the
isometries become isometries modulo homotheties.\\

Let $(K,g_K)$ be a $k$-dimensional Riemannian manifold isometrically
embedded in both the $n$-dimensional Riemannian manifolds
$(M_1,g_1)$ and $(M_2,g_2)$,
\[
\iota_i : K \hookrightarrow M_i
\]
We assume that the isometric map $\iota_1^{-1}\circ\iota_2 :
\iota_1(K) \rightarrow \iota_2(K) $ extends to a diffeomorphism
between the normal bundles of $\iota_i (K)$ in $(M_i,g_i)$, $i=1,2$.
We further assume that both the metrics $g_1$ and $g_2$ have zero
constant scalar curvature. In this section our aim is to perform a
generalized connected sum of $(M_1,g_1)$ and $(M_2,g_2)$ along
$(K,g_K)$ and to construct on the new manifold $M = M_1
\,\sharp_{K}\, M_2$ a family of metrics
$(g_{\ep})_{\ep \in (0,1)}$, whose scalar curvature is close to zero in a suitable sense.\\

For a fixed $\ep \in (0,1)$, we describe the generalized connected
sum construction and the definition of the metric $g_\ep$ in local
coordinates, the fact that this construction yields a globally
defined metric will follow at once.\\

Let $U^k$ be an open set of $\mathbb{R}^k $, $B^{m-k}$ the
$(m-k)$-dimensional open ball ($m-k \geq 3$). For $i=1,2$, $F_i :
U^k \times B^{m-k} \rightarrow W_i \subset M_i$ given by
\[ F_i(z,x) : = \exp^{M_i}_z(x)
\]
defines local Fermi coordinates near the coordinate patches
$F_i(\cdot,0)\left( U \right)\subset \iota_i(K) \subset M_i$. In
these coordinates, the metric $g_i$ can be decomposed as
\begin{eqnarray*}
g_i(z,x) & = & g^{(i)}_{j_1j_2} dz^{j_1} \otimes dz^{j_2} +
g^{(i)}_{\al \be}dx^{\al} \otimes dx^{\be} + g^{(i)}_{j \al }dz^{j}
\otimes dx^{\al}
\end{eqnarray*}
and it is well known that in this coordinate system
\begin{eqnarray*}
g^{(i)}_{\al \be} = \delta_{\al \be}+ \bigo{|x|^2} & \qquad
\mbox{and} \qquad & g^{(i)}_{j \al } = \bigo{|x|}
\end{eqnarray*}\\

In order to perform the identification between $W_1$ and $W_2$ and
in order to glue the metrics together and define $g_\ep$, we
partially change the coordinate system, by setting \[ x = \ep \,
\ex^{-t} \, \theta \] on $F_1^{-1}(W_1)$ and \[ x = \ep \, \ex^{t}
\theta \] on $F_2^{-1}(W_2)$, for $\ep \in (0,1)$, $\log\ep < t < -
\log\ep$, $\theta \in S^{m-k-1}$. \\

Using these changes of coordinates the expressions of the two
metrics $g_1$ and $g_2$ on $U^k \times A^{1}_{\ep^2}$, where $
A^{1}_{\ep^2}$ is the annulus $\{\ep^2 < |x| < 1 \}$ become
respectively
\begin{eqnarray*}
g_{1}(z,t,\te) & = & g^{(1)}_{ij} dz^i \otimes dz^j \\
& + & u_{\ep}^{(1)} \, ^{\frac{4}{n-2}} \left[ \left(dt \otimes dt +
g^{(1)}_{\lambda \mu} d\theta^\lambda \otimes d\theta^\mu \right) +
g^{(1)}_{t \theta}dt \ltimes d\theta
\right]\\
& + & g^{(1)}_{i t} d\zi \otimes dt + g^{(1)}_{i \lambda} d\zi
\otimes d\te^{\lambda}
\end{eqnarray*}
and
\begin{eqnarray*}
g_{2}(z,t,\te) & = & g^{(2)}_{ij} dz^i \otimes dz^j \\
& + & u_{\ep}^{(2)} \, ^{\frac{4}{n-2}} \left[ \left(dt \otimes dt +
g^{(2)}_{\lambda \mu} d\theta^\lambda \otimes d\theta^\mu \right) +
g^{(2)}_{t \theta}dt \ltimes d\theta
\right]\\
& + & g^{(2)}_{i t} d\zi \otimes dt + g^{(2)}_{i \lambda} d\zi
\otimes d\te^{\lambda}
\end{eqnarray*}
where by the compact notation $g_{t \theta} \, dt \ltimes d\theta$
we indicate the general component of the normal metric tensor (that
is, it involves $dt \otimes dt$, $d\theta^\lambda \otimes
d\theta^\mu$
and $dt \otimes d\theta^\lambda$ components).\\

Remark that for $j = 1,2$ we have
\[
\begin{array}{rlllllll}
g^{(j)}_{\lambda \mu} & = & \bigo{1} & \qquad \qquad &
g^{(j)}_{t \theta} & =& \bigo{|x|^2}\\[3mm]
%g^{(j)}_{ \lambda \mu}& = &\bigo{|x|^2}\\
g^{(j)}_{i t} & = & \bigo{|x|^2} & \qquad \qquad &
g^{(j)}_{i \lambda} & = & \bigo{|x|^2}\\
\end{array}
\]
and
\begin{eqnarray*}
u_{\ep}^{(1)} (t) = \ep^{\frac{n-2}{2}} \ex^{-\frac{n-2}{2}t} &
\qquad \mbox{and} \qquad & u_{\ep}^{(2)} (t) = \ep^{\frac{n-2}{2}}
\ex^{\frac{n-2}{2}t}
\end{eqnarray*}

We choose a cut-off function $\zeta : (\log\ep, -\log\ep)
\rightarrow [0,1]$ to be a non increasing smooth function which is
identically equal to $1$ in $(\log\ep, -1]$ and $0$ in
$[1,-\log\ep)$ and we choose another cut-off function $\eta :
(\log\ep, -\log\ep) \rightarrow [0,1]$ to be a non increasing smooth
function which is identically equal to $1$ in $(\log\ep, -\log\ep
-1]$ and which satisfies $\lim_{t\rightarrow -\log\ep} \eta = 0$.
Using these two cut-off functions, we can define a new normal
conformal factor $u_{\ep}$ by
\begin{eqnarray*}
u_{\ep} (t) : = \eta(t) \, u_{\ep}^{(1)} (t) + \eta(-t) \,
u_{\ep}^{(2)}(t)
\end{eqnarray*}
and the metric $g_{\ep}
%_{\ep \in (0,1)}
$ by
\begin{eqnarray}
g_{\ep}(z,t,\te) & := & \left(\zeta g^{(1)}_{ij} + (1-\zeta) g^{(2)}_{ij} \right)dz^i \otimes dz^j \nonumber \\
& + & u_{\ep}^{\frac{4}{n-2}} \left[ dt \otimes dt + \left(\zeta
g^{(1)}_{\lambda \mu} + (1-\zeta) g^{(2)}_{\lambda \mu}
\right)d\theta^\lambda \otimes d\theta^\mu \right. \nonumber \\
& & \qquad \qquad \qquad + \left. \left(\zeta g^{(1)}_{t \theta}
+ (1-\zeta) g^{(2)}_{t \theta} \right) dt \ltimes d\theta \right]\\
& + & \left( \zeta g^{(1)}_{i t} + (1-\zeta) g^{(2)}_{i t} \right)
d\zi
\otimes dt \nonumber \\
& + & \left(\zeta g^{(1)}_{i \lambda} + (1-\zeta) g^{(2)}_{i
\lambda}\right) d\zi \otimes d\te^{\lambda} \nonumber\\ \nonumber
\end{eqnarray}

Closer inspection of this expression shows that the metric $g_{\ep}$
- whose definition can be obviously completed by setting $g_{\ep}
\equiv g_1$ and $g_{\ep} \equiv g_2$ out of the "polyneck" - is a
Riemannian metric which is globally defined on the manifold $M$.\\

Following \cite{Mazz} it is immediate to obtain the estimate for the
scalar curvature of the approximate solution metric.
\begin{prop}
\label{SCALAR CURVATURE ESTIMATE} There exists a constant $C >0$
independent of $\ep \in (0,1)$ such that
\begin{eqnarray}
\label{scalar curvature estimate} |S_{g_{\ep}} | \leq C \, \ep^{-1}
\, (\cip \, t)^{1-n}
\end{eqnarray}
for $|t| \leq |\log \ep | -1$.
\end{prop}
Of course, when we consider $Rg_1$ and $Qg_2$ as initial metrics the
estimate \ref{scalar curvature estimate} still remains true, but the
constant $C$ depends now on the factors $P$ and $Q$.\\

Another useful tool we can obtain from \cite{Mazz} is the expression
for the $g_\ep$-laplacian on the polyneck
\begin{eqnarray*}
\label{locallapl} \Delta_{\gep} & = & \uep^{-\frac{4}{n-2}} \left[
\partial_t^2 + (n-2) \, \tip \left(\frac{n-2}{2}t \right)
\partial_t +  \Delta^{(\theta)}_{S^{n-1}} + \uep^{\frac{4}{n-2}}\Delta^{(z)}_{K}
 + \bigo{|x|}\Phi(\nabla, \nabla^2)\right]
\end{eqnarray*}
where $\Phi(\nabla, \nabla^2)$ is a nonlinear differential operator
involving first order and second order partial derivatives with
respect to $t$, $\theta^\lambda$ and $z^j$ and whose coefficients
are bounded uniformly on the "polyneck", as $\ep \in (0,1)$.

%The aim of the geometric construction is to perform the generalized
%connected sum $M$ of the manifolds $(M_1, g_1)$ and $(M_2, g_2)$
%along the common submanifold $(K, g_K)$ and to produce a family $\{
%g_\ep \}_\ep$ of approximate solution metrics on $M$.

\section{Analysis of a linear operator}

Our aim is now to solve the Yamabe equation
\begin{eqnarray}
\label{Yamabe} \Delta_{\gep}u + c_m S_{\gep} u & = & c_m S
u^{\frac{m+2}{m-2}}
\end{eqnarray}
where $c_m = - (m-2) / 4 (m-1)$ and $S = S(\ep)$ is a suitable
constant.\\

Since we want to preserve the structure of the two initial metrics
far away from the gluing locus, we are looking for a conformal
factor $u$ as close to $1$ as we want. For these reasons it is
natural to consider the change $u = 1 + v$ and consequently the
equation
\begin{eqnarray}
\label{Yamabe in 1} \Delta_{\gep} v   &
= & c_m \left( S- S_{\gep} \right) + c_m \left( S - S_{\gep} \right) v + c_m\frac{4}{m-2} S v +  c_m S f(v)  \\
& =: & F_\ep(v) \nonumber
\end{eqnarray}
where $ f(v) = \left((1 + v)^{\frac{m+2}{m-2}} - 1 -\frac{m+2}{m-2}v
\right)$.\\

Since the first eigenvalue $\mu_\ep$ of the operator $\Delta_{\gep}$
is a $\bigo{\ep^{n-2}}$, it is not easy to provide a good estimate
for the inverse of the laplacian dealing directly with the equation
above. It is better, on the other hand, to consider the following
problem
\begin{eqnarray}
\label{Yamabe projection} \Delta_{g_\ep} v & = & F_\ep(v) -
\lambda(\ep,v)\beta_\ep
\end{eqnarray}
where $\beta_\ep = c_1 \chi_1 - c_2 \chi_2$, with $c_1,c_2 > 0$ and
$\chi_1, \chi_2 $ smooth monotone cut-off defined by
\[
\chi_1 = \left\{
           \begin{array}{ll}
             1 \quad & \hbox{on $M_1 \setminus T^\ep_0$} \\
             1 \quad & \hbox{on $\{ \log \ep < t < \log\ep + \al  \}$} \\
             0 \quad & \hbox{on $\{ \log \ep + \al + 1 < t < 0 \}$} \\
             0 \quad & \hbox{otherwise}
           \end{array}
         \right.
\]
\[
\chi_2 = \left\{
           \begin{array}{ll}
             1 \quad & \hbox{on $M_2 \setminus T^\ep_0$} \\
             1 \quad & \hbox{on $\{ - \log \ep + \al < t < - \log\ep  \}$} \\
             0 \quad & \hbox{on $\{  0< t < -\log \ep - \al - 1 \}$} \\
             0 \quad & \hbox{otherwise}
           \end{array}
         \right.
\]
is such that $\int_M \beta_\ep \dvol{\gep} = 0$ and $\int_M
{\beta_\ep}^2 \dvol{\gep}$ and we can think of it as an
approximation of the first eigenvector of $\Delta_{\gep}$. In this
problem we are looking for a function $v$ and an approximate first
eigenvalue $\lambda(\ep,v)$ such that the equation \ref{Yamabe
projection} is verified.\\

Once this problem will be solved, we will show that, by scaling the
initial metrics $g_1$ and $g_2$, the constant $\lambda(\ep,v)$ can
be chosen to be zero, providing a solution of the equation
\ref{Yamabe in 1}.\\

By linearizing the equation \ref{Yamabe projection} we are induced
to consider the linear problem
\begin{eqnarray}
\label{linear projection}
\Delta_{\gep} u & = & f - \lambda
\beta_\ep
\end{eqnarray}
where $f$ is an assigned function such that $\int_M f \dvol{\gep}
=0$ and we are looking for a suitable constant $\lambda$ and a
solution $u$ which, up to a constant, can be chosen such that
$\int_M u \dvol{\gep} = 0$.\\

In order to choose a good functional setting for this linear
problem, let us recall the following result from \cite{Mazz}
\begin{prop}
\label{LAPE} Given $ \gamma \in (0, n-2)$, there exist a real number
$\al = \al (n, \gamma) > 0$ and a constant $C_{n, \gamma} \geq 0$
such that for all $\ep \in (0, \ex^{-\al})$ and all $v,f \in
{\mathcal C}^{0}(T^{\ep}_{\al})$ satisfying $\Delta_{\gep} v = f$,
the following estimate holds
\begin{eqnarray}
\label{lape} \sup_{T^{\ep}_{\al}} \left| \psi_{\ep}^{\gamma} v
\right| & \leq & C_{n, \gamma} \left(\,\, \sup_{T^{\ep}_{\al}}
\left| \psi_{\ep}^{\gamma + 2} f \right|\,\, + \,\, \sup_{\partial
T^{\ep}_{\al}} \left| \psi_{\ep}^{\gamma} v \right| \,\,\right)
\end{eqnarray}
where $T^{\ep}_{\rho} : = \{\log \ep + \rho \leq t \leq - \log \ep
-\rho \}$, for $\rho > 0$ and the weight $\psi_\ep$ interpolate
smoothly between these definitions in $T^{\ep}_0 \setminus
T^{\ep}_{\al}$
\[ \psi_{\ep} : =
\begin{cases}
\ep \cip t & \text{in $T^{\ep}_{\al}$} \\
1 & \text{in $M \setminus T^{\ep}_0$}
\end{cases}
\]
\end{prop}

Having this result it becomes quite natural to consider functions
$f\in \mathcal{C}^0(M)$ such that $\nor{f}{\mathcal{C}^0_{\gamma +
2}(M)} = \sup_M \left| \psi_\ep^{\gamma +2} f \right| < + \infty$
and looking for solutions $u\in \mathcal{C}^0(M)$ such that
$\nor{u}{\mathcal{C}^0_{\gamma}(M)} = \sup_M \left|
\psi_\ep^{\gamma} u \right| < + \infty$.\\

As a first step towards the solution of the problem \ref{linear
projection} we will proof the following

\begin{lemma}
\label{ERROR ESTIMATE} Given a function $f \in \mathcal{C}^0(M)$
such that $\int_M f \dvol{\gep} =0$ and $\norgd{f}{M} < + \infty$ it
is possible to find a real number $\lambda$, an approximate solution
$u \in \mathcal{C}^0(M)$ such that $\int_M u \dvol{\gep} =0$ and
$\norg{u}{M} < + \infty$ and an error term $R \in \mathcal{C}^0(M)$
such that $\int_M R \dvol{\gep} =0$ and $\norgd{R}{M} < + \infty$
that verify
\begin{eqnarray*}
\Delta_{\gep} u & = & f - \lambda\beta_\ep + R
\end{eqnarray*}
Moreover there exist positive constants $A,B,C>0$ such that the
following estimates yield for every $\gamma \in (0,n-2)$
%\begin{enumerate}
%  \item $\norgd{R}{M} \leq A \ep^{\frac{n-2}{2}} \norgd{f}{M}$
%  \item $\norg{u}{M} \leq B \norgd{f}{M}$
%  \item $|\lambda| \leq C \norgd{f}{M}$
%\end{enumerate}
\begin{eqnarray}
\label{error estimate}
\norgd{R}{M} & \leq  & A \, \ex^{(n-2)\al} \, \ep^{n-2} \norgd{f}{M} \\
\label{approx solution estimate}
\norg{u}{M} & \leq  & B \, \norgd{f}{M} \\
\label{lambda estimate} |\lambda| & \leq  & C \, \norgd{f}{M}
\end{eqnarray}
\end{lemma}

The proof of the lemma consists in building an approximate solution
$u$ and in estimating the error term. In order to do that let us
consider a non negative smooth function $\chi_P$ such that the
triple $\{ \chi_1, \chi_P, \chi_2 \}$ is a partition of the unity.
We can write
\[
f = f\chi_1 + f\chi_P + f\chi_2 = f_1 + f_P + f_2
\]
As a first step we want to build a good approximate solution on the
polyneck. It is well known that the problem
\[
\left\{
  \begin{array}{ll}
    \Delta _{\gep} v = f_P & \hbox{on $T^{\ep}_\al $ } \\
     \\
    v = 0 & \hbox{on $\partial T^{\ep}_\al $}
  \end{array}
\right.
\]
admits a solution and we call it $\tilde{u}_P$. Moreover, if $f_P$
is continuous, so does $\tilde{u}_P$ and thanks to the Lemma
\ref{LAPE}, if we choose $\al$ large enough, we have that the
following estimate yields
\begin{eqnarray*}
\norg{\tilde{u}_P}{T^{\ep}_\al} \leq C_P \, \norgd{f_P}{T^{\ep}_\al}
\end{eqnarray*}
for some positive constant $C_P>0$. Notice that the boundary
condition allows us to drop out the term
$\norg{\tilde{u}_P}{\partial T^{\ep}_\al}$ in the above estimate.\\

Let us define $u_P := \chi_P \tilde{u}_P$, as a consequence we have
that
\begin{eqnarray*}
\Delta_{\gep} u_p & = & \Delta_{\gep} \tilde{u}_P - \Delta_{\gep}
(1-\chi_P)\tilde{u}_P \\
& = & f_P - \Delta_{\gep} (\chi_1 \tup) - \Delta_{\gep} (\chi_2 \tup)\\
& = & f_P - q_1 - q_2 \\
\end{eqnarray*}
where $q_i = \Delta_{\gep} (\chi_i \tup)$, $i = 1,2$.\\

Let us call $\tf_i := f_i + q_i$, $i=1,2$ and $\tf := \tf_1 +
\tf_2$. Since $\int_M f \dvol{\gep} = 0$, it is easy to check that
also $\int_M \tf \dvol{\gep} = 0$. Hence $\int_M \tf_1 \dvol{\gep} =
- \int_M \tf_2 \dvol{\gep}$.\\

Let us define $h_i := \tf_i + (-1)^i \lambda c_i \chi_i$ for $i=1,2$
and $h := h_1 + h_2 = \tf - \lambda \beta_\ep$. Obviously we have
that $\int_M h \dvol{\gep} = 0 $ and $\int_M h_1 \dvol{\gep} = -
\int_M h_2 \dvol{\gep}$.\\

Moreover
$$\int_M h_1 \dvol{\gep} - \int_M h_2 \dvol{\gep} =  \int_M \tf_1 - \tf_2 \dvol{\gep}
- \lambda \int_M c_1 \chi_1 + c_2 \chi_2 \dvol{\gep}$$ Hence, by
setting $$\lambda := \frac{\int_M \tf_1 - \tf_2 \dvol{\gep}}{\int_M
c_1 \chi_1 + c_2 \chi_2 \dvol{\gep}}$$ it follows at once that
$\int_M h_i \dvol{\gep} = 0$, $i=1,2$.\\

As a second step we want now to construct approximate solutions on
the pieces of $M$ coming from $M_1$ and $M_2$. For this purpose, let
us consider, for $i=1,2$, the functions $\tu_i$ verifying
\begin{eqnarray*}
\Delta_{g_i} \tu_i = h_i - b_i \delta_K
\end{eqnarray*}
where $b_i = C_{n,K} \int_M h_i \dvol{g_i}$ and $C_{n,K}$ is a
suitable constant.\\

It is rather simple to describe how this functions approximately
look like, in fact we can write (notice that the following remarks
still hold for $i=2$)
\begin{eqnarray*}
\Delta_{g_1} \tu_1 = h_1 + \frac{1}{V_1}\int_M h_1 \dvol{g_1} -
\frac{1}{V_1}\int_M h_1 \dvol{g_1} - b_1 \delta_K
\end{eqnarray*}
we can now consider the split $\tu_1 = \bar{u}_1 + \hat{u}_1$ where
\[
\left\{
  \begin{array}{ll}
    \Delta_{g_1} \bar{u}_1 = h_1 - \frac{1}{V_1}\int_M h_1 \dvol{g_1}\\
      \\
    \Delta_{g_1} \hat{u}_1 = \frac{1}{V_1}\int_M h_1 \dvol{g_1} - b_1
\delta_K
  \end{array}
\right.
\]
we can think of $\bar{u}_1$ as the "finite part" and of $\hat{u}_1$
as the "pure Green function part" of $\tu_1$. In particular
$\hat{u}_1$ has the following shape:
\begin{eqnarray*}
\hat{u}_1 & = & \Omega_{n,K}\int_M h_1 \dvol{g_1} \left( |x|^{2-n} +
\bigo{|x|^{1-n}} \right)\\
& = & \Omega_{n,K}\int_M h_1
\left(\dvol{g_1}-\dvol{\gep}\right)\left( |x|^{2-n} +
\bigo{|x|^{1-n}} \right) \\
\end{eqnarray*}
since $ \int_M h_1 \dvol{\gep} = 0 $. For the gradient $\nabla
\hat{u}_1$ the expression is
\begin{eqnarray*}
\nabla \hat{u}_1 & = & \Omega_{n,K}\int_M h_1
\left(\dvol{g_1}-\dvol{\gep}\right)\left( (n-2)|x|^{2-n} +
\bigo{|x|^{1-n}} \right)
\end{eqnarray*}
because $\nabla |x|^{2-n} = \nabla (\ep^{2-n} \ex^{(n-2)t}) =
\ep^{2-n} (n-2) \ex^{(n-2)t} = (n-2) |x|^{2-n}$.\\

In order to glue together $\tu_1$ and $\tu_2$, we will use the
following smooth and monotone cut-off functions
\[
\phi_1 = \left\{
           \begin{array}{ll}
             1 \quad & \hbox{on $M_1 \setminus T^\ep_0$} \\
             1 \quad & \hbox{on $\{ \log \ep < t < \log\ep + \al + 1  \}$} \\
             0 \quad & \hbox{on $\{ \log \ep + \al + 2 < t < 0 \}$} \\
             0 \quad & \hbox{otherwise}
           \end{array}
         \right.
\]
\[
\phi_2 = \left\{
           \begin{array}{ll}
             1 \quad & \hbox{on $M_2 \setminus T^\ep_0$} \\
             1 \quad & \hbox{on $\{ - \log \ep - \al - 1 < t < -\log\ep  \}$} \\
             0 \quad & \hbox{on $\{  0< t < -\log \ep -\al - 2 \}$} \\
             0 \quad & \hbox{otherwise}
           \end{array}
         \right.
\]
Now we can define the approximate solution $u$ as
$$ u = u_P + \phi_1 \tu_1 + \phi_2 \tu_2 $$
and we can calculate
\begin{eqnarray*}
\Delta_{\gep} u & = & \Delta_{\gep} u_P + \Delta_{\gep}(\phi_1
\tu_1) + \Delta_{\gep} (\phi_2 \tu_2) \\
& = & f_P - q_1 - q_2 \\
& + & \phi_1 (h_1 - b_1 \delta_K) + (\Delta_{\gep} \phi_1) \tu_1 +
\gep (\nabla \phi_1 , \nabla \tu_1) \\
& + & \phi_2 (h_2 - b_2 \delta_K) + (\Delta_{\gep} \phi_2) \tu_2 +
\gep (\nabla \phi_2 , \nabla \tu_2)\\
& = & f - \lambda \beta_\ep\\
& + & (\Delta_{\gep} \phi_1) \tu_1 +
\gep (\nabla \phi_1 , \nabla \tu_1) \\
& + & (\Delta_{\gep} \phi_2) \tu_2 +
\gep (\nabla \phi_2 , \nabla \tu_2)\\
\end{eqnarray*}
At this point it is quite natural to define $E_i := (\Delta_{\gep}
\phi_i) \tu_i + \gep (\nabla \phi_i , \nabla \tu_i)$, $i=1,2$  and
$R := E_1 + E_2$, so that
\begin{eqnarray*}
\Delta_{\gep} u & =& f - \lambda \beta_\ep + R
\end{eqnarray*}
We can now proceed with the estimate of $R$. Without loss of
generality, let us look for example at the error term $E_1$. Since
$\textrm{supp} (\Delta_{\gep}\phi_1)$ and $\textrm{supp} (\nabla
\phi_1) $ are both included in $[\log\ep + \al + 1, \log\ep + \al +
2]$, the term $E_1$ is supported here as well.\\

It follows from a straightforward computation that
\begin{eqnarray*}
\left| \int_M h_1 (\dvol{\gep} - \dvol{g_1})\right| & \leq & C^{(0)}
\, \ex^{\al\gamma} \, \ep^{n-2} \, \norgd{f}{M}
\end{eqnarray*}
for some constant $C^{(0)}>0$.\\

According to the splitting of $\tu_1$ we have that
$$ \tu_1 = \bar{u}_1 + \hat{u}_1 \qquad \qquad \nabla\tu_1 = \nabla\bar{u}_1 + \nabla \hat{u}_1$$
and by remembering the expressions found for $\hat{u}_1$ and $\nabla
\hat{u}_1$ it is easy to see that both these terms are
$\bigo{\ex^{(n-2+\gamma)\al}}$ on the support of $E_1$, hence, for
sufficiently large $\al > 0$ it is clear that $|\tu_1| < A^{(0)}
|\hat{u}_1|$ and $|\nabla \tu_1| < A^{(1)} |\nabla \hat{u}_1|$, for
suitable constants $A^{(0)}, A^{(1)}>0$.\\

Now we are ready to estimate on the interval $[\log\ep + \al + 1,
\log\ep + \al + 2]$ the term
\begin{eqnarray*}
|\psi_\ep^{\gamma + 2} (\Delta_{\gep} \phi_1) \tu_1 | & \leq &
A^{(0)}
|\psi_\ep^{\gamma + 2} (\Delta_{\gep} \phi_1) \hat{u}_1 |\\
& \leq & A^{(2)} |x|^{\gamma +2} |x|^{-2} |\Omega_{n,K}| \left|
\int_M
h_1 (\dvol{\gep} - \dvol{g_1})\right| |x|^{2-n}\\
& \leq & A^{(3)} |x|^{-(n-2-\gamma)} \ex^{\al\gamma} \ep^{n-2} \norgd{f}{M}\\
&\leq & A^{(4)} \ex^{(n-2)\al} \ep^{n-2} \norgd{f}{M}
\end{eqnarray*}

Analogously
\begin{eqnarray*}
\left|\psi_\ep^{\gamma + 2} {\gep} \left(\nabla \phi_1, \nabla \tu_1
\right)\right| & \leq & A^{(5)} \ex^{(n-2)\al} \ep^{n-2}
\norgd{f}{M}
\end{eqnarray*}

Hence
\begin{eqnarray*}
\norgd{E_1}{M} & \leq & A^{(6)} \ex^{(n-2)\al} \ep^{n-2} \norgd{f}{M}\\
\end{eqnarray*}
Since the same estimate holds for $E_2$, we conclude that, for a
suitable constant $A>0$
\begin{eqnarray*}
\norgd{R}{M} & \leq & A \ex^{(n-2)\al} \ep^{n-2} \norgd{f}{M}
\end{eqnarray*}

In order to obtain the estimate \ref{approx solution estimate}, let
us recall that $u = u_P + \phi_1 \tu_1 + \phi_2 \tu_2$ and that, for
$\al > 0$ large enough, $\norg{u_p}{T^\ep_\al} = \norg{\chi_P
\tu_P}{T^\ep_\al} \leq B^{(0)} \norgd{f_p}{T^\ep_\al}$.\\

On the other hand, on the support of $\phi_1$ we have that
\begin{eqnarray*}
|\psi_\ep^\gamma \phi_1 \tu_1| & \leq & B^{(1)} |x|^{\gamma}
|\Omega_{n,K}|
\left| \int_M h_1 \dvol{\gep} \right| |x|^{2-n} \\
& \leq & B^{(2)} \ex^{\al\gamma} \ep^{n-2} \norgd{f}{M}
|x|^{2 - n + \gamma} \\
& \leq & B^{(3)} \ex^{(n-2)\al} \ep^{\gamma} \norgd{f}{M}
\end{eqnarray*}

Hence, it is clear that there exists a constant $B>0$ such that:
\begin{eqnarray*}
\norg{u}{M} & \leq & B \norgd{f}{M}
\end{eqnarray*}

Finally, by remembering the expression for $\lambda$ it follows from
a straightforward computation that, for large enough $C>0$,
\begin{eqnarray*}
|\lambda| & \leq & C \norgd{f}{M}
\end{eqnarray*}
and the lemma \ref{ERROR ESTIMATE} is proved.\\

The idea is to solve the equation \ref{linear projection} by means
of a sequence method. We start by setting $f^{(0)} := f$ and thanks
to the lemma \ref{ERROR ESTIMATE} we obtain a triple
$(\lambda^{(0)}, u^{(0)}, R^{(0)})$ verifying the equation
$$\Delta_{\gep} u^{(0)} = f^{(0)} - \lambda^{(0)} \beta_\ep +
R^{(0)}$$ and the estimates \ref{error estimate}, \ref{approx
solution estimate}, \ref{lambda estimate}. Now we set $f^{(1)}:=
-R^{(0)}$ and we find another triple $(\lambda^{(1)}, u^{(1)},
R^{(1)})$ with the same properties of the first one and so on. In
general, for every $j \in \mathbb{N}$, we have $f^{(j)} :=
-R^{(j-1)}$ and a triple $(\lambda^{(j)}, u^{(j)}, R^{(j)})$
verifying the equation
\begin{eqnarray}
\label{j-equation} \Delta_{\gep} u^{(j)} & = & f^{(j)} -
\lambda^{(j)} \beta_\ep + R^{(j)}
\end{eqnarray}
and the estimates \ref{error estimate}, \ref{approx solution
estimate}, \ref{lambda estimate}.\\

By taking the sum of the equations \ref{j-equation} we have that,
for every $ N \in \mathbb{N} $
\begin{eqnarray*}
\Delta_{\gep} \left( \sum_{j=0}^N u^{(j)}  \right)  & = &
\sum_{j=0}^N f^{(j)} - \left( \sum_{j=0}^N  \lambda^{(j)}\right)
\beta_\ep +
\sum_{j=0}^N R^{(j)}\\
& = & f - \left( \sum_{j=0}^N  \lambda^{(j)}\right) \beta_\ep  +
R^{(N)}
\end{eqnarray*}

In other words
\begin{eqnarray*}
\Delta_{\gep} v^{(N)} & =& f - \mu^{(N)} \beta_\ep + R^{(N)}
\end{eqnarray*}
where $v^{(N)} := \sum_{j=0}^N u^{(j)} $ and $\mu^{(N)} :=
\sum_{j=0}^N  \lambda^{(j)}$.\\

Notice that from the estimate of the lemma \ref{ERROR ESTIMATE} it
follows easily that
\begin{eqnarray*}
\norgd{f^{(j)}}{M}\, = \,\norgd{R^{(j-1)}}{M} & \leq & A
\ex^{(n-2)\al}\ep^{n-2} \norgd{f^{(j-1)}}{M} \\
& \leq & \left(A \ex^{(n-2)\al} \ep^{n-2} \right)^j \norgd{f}{M} \\
& & \\
\norg{u^{(j)}}{M} & \leq & B \norgd{f^{(j)}}{M}\\
& \leq & B \left(A \ex^{(n-2)\al} \ep^{n-2} \right)^j \norgd{f}{M}\\
& & \\
|\lambda^{(j)}| & \leq & C \norgd{f^{(j)}}{M} \\
& \leq & C \left(A \ex^{(n-2)\al} \ep^{n-2} \right)^j \norgd{f}{M}
\end{eqnarray*}

Now it is clear that, for sufficiently small $\ep>0$, there exist
$\lambda \in \mathbb{R}$ and a continuous function $u$ such that
\begin{eqnarray*}
R^{(N)} & \stackrel{\nogd{\cdot}}{\longrightarrow} &  0\\
& & \\
v^{(N)} & \stackrel{\nog{\cdot}}{\longrightarrow} & u \\
& & \\
\mu^{(N)} & \longrightarrow & \lambda
\end{eqnarray*}
Moreover there exist positive constants $B',C'>0$ such that
$$\norg{u}{M}\,\, \leq \,\, B' \norgd{f}{M}  \qquad\qquad  |\lambda| \,\, \leq  \,\, C' \norgd{f}{M}$$
Hence
\[
\left\{
  \begin{array}{ll}
     \Delta_{\gep} v^{(N)} \,\, \stackrel{\nogd{\cdot}}{\longrightarrow}  \,\,  f - \lambda \beta_\ep &  \\
         &  \\
     v^{(N)} \,\,  \stackrel{\nog{\cdot}}{\longrightarrow} \,\,  u  &
  \end{array}
\right.
\]
On the other hand we have that, for every $N \in \mathbb{N}$ and for
every $\phi \in \mathcal{C}^{\infty}(M)$
\begin{eqnarray*}
\int_M v^{(N)} \Delta_{\gep} \phi \dvol{\gep} & = & \int_M (f -
\mu^{(N)}\beta_\ep + R^{(N)}) \phi \dvol{\gep}
\end{eqnarray*}
Hence, by taking the limit for $N \rightarrow + \infty$ we find the
expression
\begin{eqnarray*}
\int_M u \Delta_{\gep} \phi \dvol{\gep} & = & \int_M (f -
\lambda\beta_\ep ) \phi \dvol{\gep}
\end{eqnarray*}
for every $\phi \in \mathcal{C}^{\infty}(M)$. It is to say that
$\Delta_{\gep} u = f - \lambda \beta_\ep$ in the sense of the
distributions.\\

Thanks to the elliptic regularity (see for example \cite{Aub},
\cite{Gil-Tru}) if we suppose that $f \in \mathcal{C}^{\infty}(M)$,
then so does $u$ and the expression above is a pointwise identity.\\

To conclude this section we summarize our results in the following
\begin{prop}
\label{LINEAR PROBLEM} Given a function $f \in \mathcal{C}^0(M)$
such that $\int_M f \dvol{\gep} = 0$ and $\norgd{f}{M} < + \infty$,
it is possible to find a real number $\lambda$ and a function $u \in
\mathcal{C}^0(M)$ with $\int_M u \dvol{\gep} = 0$ and $\norg{u}{M} <
+ \infty$ verifying
\begin{eqnarray*}
\Delta_{\gep} u & = & f - \lambda \beta_\ep
\end{eqnarray*}
in the sense of the distributions and the following estimates
\begin{eqnarray}
\label{linear solution estimate}
\norg{u}{M} & \leq & B' \norgd{f}{M} \\
|\lambda| & \leq & C' \norgd{f}{M}
\end{eqnarray}
for suitable constants $B',C'>0$.\\

Moreover, if $f \in \mathcal{C}^{\infty}(M)$, then $u \in
\mathcal{C}^{\infty}(M)$ and the identity above holds pointwise.
\end{prop}

\section{Fixed point argument}

The aim of this section is to solve the problem \ref{Yamabe
projection}
\begin{eqnarray*}
\Delta_{g_\ep} v & = & F_\ep(v) - \lambda(\ep,v)\beta_\ep
\end{eqnarray*}
We will be able to do that by means of the results of the previous
section and of a contracting mapping argument.\\

Before starting, let us remark that in the expression for $F_\ep
(v)$ (see equation \ref{Yamabe in 1}) it is always possible to
choose $S = S(\ep,v)$ in such a way that $\int_M F_\ep \dvol{\gep} =
0$. Moreover, by using the scalar curvature estimate $S_{\gep} =
\bigo{\ep^{n-2} |x|^{1-n}}$ found in \cite{Mazz} it will be easy to
see that $S = \bigo{\ep^{n-2}}$.\\

Now, for every $\sigma>0$ let us define the space
$$\mathcal{C}^0_\sigma (M) \,\, := \,\, \left\{ h \in \mathcal{C}^0 (M) \,\,
: \,\, \int_M h \dvol{\gep} = 0 \quad \hbox{and} \quad
\nor{h}{\mathcal{C}^0_\sigma (M)} < + \infty  \right\}$$ where
$\nor{h}{\mathcal{C}^0_\sigma (M)} := \sup_M \left|\psi_\ep^{\sigma}
h \right|$.\\

Let us define also the maps
\begin{eqnarray*}
H_\ep \,\, : \,\, \mathcal{C}^0_\gamma (M) & \longrightarrow &
\mathcal{C}^0_{\gamma + 2} (M) \\
v & \longmapsto & F_\ep (v) - \lambda (\ep,v) \beta_\ep\\
& & \\
\Delta_{\gep}^{-1} \,\, : \,\, \mathcal{C}^0_{\gamma + 2} (M) &
\longrightarrow & \mathcal{C}^0_\gamma (M)\\
w & \longmapsto & \Delta_{\gep}^{-1} w \\
& & \\
P_\ep \,\, : \,\, \mathcal{C}^0_\gamma (M) & \longrightarrow &
\mathcal{C}^0_\gamma(M)\\
v & \longmapsto & \Delta_{\gep}^{-1} \circ H_\ep (v)
\end{eqnarray*}

Let us start with the following lemma
\begin{lemma}
\label{Ball in Ball} For $\gamma \in (0,1)$ and for sufficiently
small $\ep>0$ there exists a radius $r_\ep>0$ such that $P_\ep
\left( B^\gamma_{r_\ep} \right)\subset B^\gamma_{r_\ep} $, where
$B^\gamma_{r_\ep}:= \left\{ u \in \mathcal{C}^0_\gamma (M) \,\, :
\,\, \norg{u}{M} \leq r_\ep \right\}$. In other words:
\begin{eqnarray}
\norg{v}{M} \leq r_\ep & \Longrightarrow & \norg{P_\ep (v)}{M} \leq
r_\ep
\end{eqnarray}
\end{lemma}
In order to prove the statement we observe that from the estimate
\ref{linear solution estimate} of the proposition \ref{LINEAR
PROBLEM} we obtain immediately the inequality
\begin{eqnarray*}
\norg{P_\ep (v)}{M} & \leq & D \norgd{F_\ep(v)}{M}
\end{eqnarray*}
for large enough $D>0$.\\

Now we have to estimate the term
\begin{eqnarray*}
|F_\ep (v) \psi^{\gamma + 2}_\ep| & \leq & c_m |S| \psi^{\gamma +
2}_\ep + c_m |S_{\gep}| \psi^{\gamma}_\ep\\
& + & c_m |S| \psi^{\gamma + 2}_\ep  |v| + c_m |S_{\gep}|
\psi^{\gamma}_\ep |v| + \frac{|S|}{m-1} \psi^{\gamma}_\ep |v| \\
& + & c_m |S| |f(v)| \psi^{\gamma}_\ep\\
& \leq & D_{(1)} \ep^{n-2} + D_{(2)} \left[ \ep^{n-2} + \ep^{1 +
\gamma} \right] \\
& + & D_{(3)} \ep^{n-2} r_\ep + D_{(4)} \left[ \ep^{n-2} + \ep
\right] r_\ep + D_{(5)} \ep^{n-2} r_\ep \\
& + & D_{(6)} \ep^{n-2} r_\ep^2 \\
& \leq & D_{(7)} \ep^{n-2} + D_{(8)} \ep^{1 + \gamma}\\
& + & D_{(9)} \ep^{n-2} r_\ep + D_{(10)} \ep r_\ep \\
& + & D_{(11)} \ep^{n-2} r_\ep^2
\end{eqnarray*}
for suitable constants $D_{(j)}$.\\

If $n=3$, then $\ep^{n-2} = \ep
>> \ep^{1 + \gamma}$, when $\ep>0$ is small. Hence $D_{(7)}
\ep^{n-2} + D_{(8)} \ep^{1 + \gamma} \leq D_{(12)} \ep$, so if we
define $r_\ep := 2DD_{(12)} \ep$ we obtain $ D_{(9)} \ep^{n-2} +
D_{(10)} \ep + D_{(11)} \ep^{n-2} r_\ep \leq D_{(13)} \ep + D_{(11)}
\ep 2DD_{(12)} \ep \leq 1 / 2D$, for small enough $\ep > 0$ and that
guarantees $D \norgd{F_\ep(v)}{M} \leq r_\ep$.\\

If $n \geq 4$, then $\ep^{n-2} << \ep^{1 + \gamma}$ (since $\gamma
\in (0,1)$), when $\ep>0$ is small. Hence $D_{(7)} \ep^{n-2} +
D_{(8)} \ep^{1 + \gamma} \leq D_{(12)} \ep^{1+ \gamma}$, so if we
define $r_\ep := 2DD_{(14)} \ep^{1+ \gamma}$ we obtain $ D_{(9)}
\ep^{n-2} + D_{(10)} \ep + D_{(11)} \ep^{n-2} r_\ep \leq D_{(15)}
\ep + D_{(11)} \ep 2DD_{(14)} \ep^{1+\gamma} \leq 1 / 2D$, for small
enough $\ep> 0$ and that still guarantees $D \norgd{F_\ep(v)}{M}
\leq r_\ep$.\\

In both the cases the lemma is proved.\\

At this point our purpose is to prove the convergence of the
sequence $v^j := P_\ep^j (0)$ with respect to the norm
$\norg{\cdot}{M}$. Towards this aim we need to provide an estimate
of $\norg{P_\ep(u) - P_\ep(v)}{M}$ in terms of $\norg{u - v}{M}$,
where $u,v \in B^\gamma_{r_\ep}$; in fact, since $0 \in
B^\gamma_{r_\ep}$, all the terms of the sequence lie in
$B^\gamma_{r_\ep}$, because of the lemma \ref{Ball in Ball}.\\

Since $\Delta_{\gep} \left( P_\ep(u) - P_\ep(v) \right) = H_\ep (u)
- H_\ep(v)$ we have immediately that
\begin{eqnarray*}
\norg{P_\ep(u) - P_\ep(v)}{M} & \leq & C' \norgd{H_\ep(u) -
H_\ep(v)}{M} \\
& = & C' \norgd{F_\ep(u) - F_\ep(v) - \left( \lambda(\ep,u) -
\lambda(\ep, v) \right) \beta_\ep }{M}
\end{eqnarray*}

On the other hand it is easy to check that the mapping $f
\longmapsto \lambda_f$, where $f$ and $\lambda = \lambda_f$ are
those of the proposition \ref{LINEAR PROBLEM}, is a linear mapping,
therefore
$$ \lambda(\ep,u) - \lambda(\ep, v) =: \lambda_{F_\ep(u)} - \lambda_{F_\ep(v)}
= \lambda_{F_\ep(u) - F_\ep(v)} $$ Hence, thanks to the estimate
\ref{lambda estimate}, we obtain
\begin{eqnarray*}
\norg{P_\ep(u) - P_\ep(v)}{M} & \leq & C'' \norgd{F_\ep(u) -
F_\ep(v)}{M}
\end{eqnarray*}

Since the function $f$ that appears in the definition of $F_\ep(v)$
verifies the following inequality
\begin{eqnarray*}
|f(u) - f(v)| & \leq & \left[ A (|u| + |v|)  +  B
\left(|u|^{\frac{4}{m-2}}  - |v| ^{\frac{4}{m-2}}\right) \right]
|u-v|
\end{eqnarray*}
for suitable constants $A,B>0$, we can proceed to the estimate of
the term $F_\ep(u) - F_\ep(v)$: $\gamma \in (0, 1/2)$ is a
sufficient condition to ensure that
\begin{eqnarray*}
\psi_\ep^{\gamma + 2}|F_\ep(u) - F_\ep(v)|  & \leq & \psi_
\ep^{\gamma + 2} \left[ c_m |S| |u-v| + c_m |S_{\gep}| |u-v| +
\frac{|S|}{m-1} |u - v| \right] \\
& + &  \psi_\ep^{\gamma + 2} c_m |S| \left[ A (|u| + |v|)  +  B
\left(|u|^{\frac{4}{m-2}} - |v| ^{\frac{4}{m-2}}\right) \right]
|u-v|\\
& \leq & C''' \ep \norg{u-v}{M}\\
& + & C'''' \ep^{n-2} \left(  \norg{u}{M} + \norg{v}{M} \right)
\norg{u-v}{M} \\
& + & C''''' \ep^{n-2} \left( \norg{u}{M}^{\frac{4}{m-2}} +
\norg{v}{M}^{\frac{4}{m-2}} \right) \norg{u-v}{M}
\end{eqnarray*}
Hence, for $u,v \in B^\gamma_{r_\ep}$ and small enough $\ep>0$ we
get the inequality
\begin{eqnarray*}
\norg{P_\ep(u) - P_\ep(v)}{M} & \leq & C \ep \norg{u - v}{M}
\end{eqnarray*}
Now, for integers $p \leq q $ we have that
\begin{eqnarray*}
\norg{v^q - v^p}{M} & \leq & \sum_1^{p-q} \norg{v^{p+j} -
v^{p+j-1}}{M} \\
& \leq & (C\ep)^p \sum_0^{+\infty} (C\ep)^j \norg{v^1 - v^0}{M}
\end{eqnarray*}

Hence the sequence $(v^j)$ is a Cauchy sequence and it must converge
to a continuous function $v_\ep \in B^\gamma_{r_\ep}$ which is the
fixed point we were looking for i. e.
\begin{eqnarray}
P_\ep (v_\ep) & = & v_\ep
\end{eqnarray}
In other words
\begin{eqnarray}
\Delta_{\gep} v_\ep & = & F_\ep (v_\ep) - \lambda_{F_\ep(v_\ep)}
\beta_\ep
\end{eqnarray}

By means of a classical boot strap argument it is easy to see that
$v_\ep$ is actually a smooth function.

\section{The approximate eigenvalue $\lambda_{F_\ep(v_\ep)}$}

In this section we want to study the sign of the approximate
eigenvalue $\lambda_{F_\ep(v_\ep)}$. In particular our purpose is to
show that by moving the initial metrics (more precisely by scaling
up or down $g_1$ and $g_2$) the approximate eigenvalue becomes
positive or negative, hence there exist suitable constants $R, Q >
0$ such that the construction starting by $R g_1$ and $Q g_2$ as
initial metrics has zero approximate eigenvalue. Therefore, in this
case, $v_\ep$ is a solution of the problem \ref{Yamabe in 1}.\\

Since an explicit expression of $\lambda_{F_\ep(v_\ep)}$ in terms of
the initial metrics is not available, we have to handle with its
approximations, taking care in estimating the errors. Thanks to the
proposition \ref{LINEAR PROBLEM} we can think of
$\lambda_{F_\ep(v_\ep)}$ as obtained by a sequence method, exactly
like the real number $\lambda$ of the mentioned proposition.
Therefore it is possible to find a sequence
$\lambda_{F_\ep(v_\ep)}^{(j)}$ such that
\begin{eqnarray*}
\lambda_{F_\ep(v_\ep)} & = & \sum_{j=0}^{\infty}
\lambda_{F_\ep(v_\ep)}^{(j)}
\end{eqnarray*}
and such that the following estimate holds:
\begin{eqnarray*}
\left| \lambda_{F_\ep(v_\ep)}^{(j)} \right| & \leq & C \left( A
\ex^{(n-2)\al} \ep^{n-2} \right)^j \norgd{F_{\ep}(v_\ep)}{M}
\end{eqnarray*}
for suitable constants $A,C > 0$.\\

Now, for sufficiently small $\ep > 0$ it is quite easy to estimate
the difference
\begin{eqnarray*}
\left|\lambda_{F_\ep(v_\ep)} - \lambda_{F_\ep(v_\ep)}^{(0)} \right|
& \leq & \sum_{j=1}^{+\infty}
\left|\lambda_{F_\ep(v_\ep)}^{(j)}\right| \\
& \leq & C \left[ \sum_{j=1}^{+\infty} \left( A \ex^{(n-2)\al}
\ep^{n-2} \right)^j \right] \norgd{F_{\ep}(v_\ep)}{M} \\
& \leq & \frac{A \ex^{(n-2)\al} \ep^{n-2}}{1 - A \ex^{(n-2)\al}
\ep^{n-2}} C \norgd{F_{\ep}(v_\ep)}{M} \\
& \leq & B  \ex^{(n-2)\al} \ep^{n-2} \norgd{F_{\ep}(v_\ep)}{M}
\end{eqnarray*}
where $B>0$ is a suitable constant.\\

At the moment we have obtained that
\begin{eqnarray*}
\left| \lambda_{F_\ep(v_\ep)} \right| & \leq &
\left|\lambda_{F_\ep(v_\ep)}^{(0)}\right| - B \ex^{(n-2)\al}
\ep^{n-2} \norgd{F_{\ep}(v_\ep)}{M}
\end{eqnarray*}

Following the proof of the proposition \ref{LINEAR PROBLEM}, we can
write down the expression for
$\left|\lambda_{F_\ep(v_\ep)}^{(0)}\right|$, it is to say:
\begin{eqnarray*}
\lambda_{F_\ep(v_\ep)}^{(0)} & = & \frac{\int_M \widetilde{F_\ep
(v_\ep)}_1  -  \widetilde{F_\ep(v_\ep)}_2  \dvol{\gep} }{ \int_M c_1
\chi_1  +  c_2 \chi_2  \dvol{\gep} }
\end{eqnarray*}
where $\widetilde{F_\ep (v_\ep)}_i = \chi_i F_\ep (v_\ep) +
\Delta_{\gep} \left( \chi_i \tup^\ep  \right)$, for $i = 1,2$ and
where $\tup^\ep$ is the solution of the problem
\[
\left\{
  \begin{array}{ll}
    \Delta_{\gep} \tup^\ep = \chi_P F_\ep(v_\ep) & \hbox{on $T^{\ep}_\al$} \\
     \\
    \tup^\ep = 0 & \hbox{on $\partial T^{\ep}_\al$}
  \end{array}
\right.
\]

It is convenient to write
\begin{eqnarray}
\label{lambdeffepsilonzero}
\lambda_{F_\ep(v_\ep)}^{(0)} & = &
\frac{1}{ \int_M c_1 \chi_1 + c_2 \chi_2  \dvol{\gep} }   \int_M
{F_\ep (v_\ep)} \chi_1 -
{F_\ep(v_\ep)}\chi_2  \dvol{\gep} \\
& + & \frac{1}{ \int_M c_1 \chi_1 + c_2 \chi_2  \dvol{\gep} } \int_M
\Delta_{\gep} (\chi_1 \tup^\ep) - \Delta_{\gep}(\chi_2 \tup^\ep)
\dvol{\gep} \nonumber
\end{eqnarray}

Concerning the first summand, it is sufficient to remember the
estimate of the scalar curvature contained in the proposition
\ref{SCALAR CURVATURE ESTIMATE} to conclude that
\begin{eqnarray*}
\int_M {F_\ep (v_\ep)} \chi_1 - {F_\ep(v_\ep)}\chi_2  \dvol{\gep} &
= & \bigo{\ep^{n-2}}
\end{eqnarray*}

Concerning the second summand, it is useful to consider, for $i=1,2$
the split $\Delta_{\gep}\left(\chi_i \tup^\ep\right) = \chi_i \left(
\Delta_{\gep} \tup^\ep \right) + \gep\left( \nabla \chi_i , \nabla
\tup^\ep \right)  +  \tup^\ep \left(  \Delta_{\gep} \chi_i \right)$.
Obviously, we have that:
\begin{eqnarray*}
\int_M \chi_i \left(\Delta_{\gep} \tup^\ep \right) \dvol{\gep} & = &
\int_{ \{\log \ep + \al \leq t \leq \log \ep + \al + 1\} } \chi_1
\chi_P F_\ep(v_\ep) \dvol{\gep} \\
& = & \bigo{ \ex^{-\al} \ep^{n-2} }
\end{eqnarray*}
and, using the Green formula, we can get now
\begin{eqnarray*}
\int_M \tup^\ep \left(  \Delta_{\gep} \chi_i \right) \dvol{\gep} +
\int_M \gep\left( \nabla \chi_i , \nabla \tup^\ep \right)
\dvol{\gep} & = & \int_{T^\ep_\al} \tup^\ep \left(  \Delta_{\gep}
\chi_i \right) \dvol{\gep} \\
& + & \int_{T^\ep_\al} \nabla \chi_i (\tup^\ep) \dvol{\gep} \\
& = & \int_{\partial T^\ep_\al} \tup^{\ep} \partial_t(\chi_i)
\dvol{\iota^* \gep} \\
& = & 0
\end{eqnarray*}
where $\iota : \partial T^\ep_\al \rightarrow M$ is the natural
embedding.\\

Hence the second summand in \ref{lambdeffepsilonzero} is a
$\bigo{\ex^{-\al} \ep^{n-2}}$.\\

Now, by performing the right choice of $\al = \al (\ep)$ we will be
able to show that the sign of the approximate eigenvalue is
determined by the sign of $\lambda_{F_\ep(v_\ep)}^{(0)}$ and, in
particular, by the term
$$ \frac{1}{ \int_M c_1 \chi_1 + c_2 \chi_2
\dvol{\gep} }   \int_M {F_\ep (v_\ep)} \chi_1 - {F_\ep(v_\ep)}\chi_2
\dvol{\gep} $$
If we set, for example $\al = -\log \ep / 2(n-2)$, we
get immediately
\begin{eqnarray*}
\ex^{(n-2)\al} \ep^{n-2} \norgd{F_\ep(v_\ep)}{M} & = & \bigo{ \ep^{n - 3/2} } \\
& & \\
\frac{1}{ \int_M c_1 \chi_1 + c_2 \chi_2  \dvol{\gep} }
\int_M \Delta_{\gep} (\chi_1 \tup^\ep) - \Delta_{\gep}(\chi_2
\tup^\ep) \dvol{\gep} & = & \bigo{ \ep^{n-2 + 1/2(n-2)} }
\end{eqnarray*}
Hence, for small $\ep > 0$ the leading term is the one we wished.\\

More precisely, if we look at the expression for $F_\ep (v_\ep)$, it
is clear that, when $\ep$ is close to zero, the sign of the
approximate eigenvalue is determined by the term
\begin{eqnarray*}
\int_M \left(S - S_{\gep} \right) \chi_1 \dvol{\gep} - \int_M (S -
S_{\gep} ) \chi_2 \dvol{\gep}
\end{eqnarray*}

At this moment we want to replace for instance the initial metric
$g_1$ by its homothetic $R g_1$ and to show that for sufficiently
large $R>0$ the sign of the expression above is determined.\\

If we indicate by $(\cdot)^R$ the geometric quantities obtained in
this case, we find the expansions
\begin{eqnarray*}
\chi_1 \dvol{\gep}^R & = & R^{\frac{m}{2}} \dvol{g_1} +
R^{\frac{k}{2} + \frac{n+2}{4} } \bigo{\ep^n \ex^{-2t}} \di t \di
\theta^1 \ldots \di \theta^{n-1}\di z^1 \ldots \di z^{k} \\
\chi_2 \dvol{\gep}^R & = & \dvol{g_2} + R^{\frac{n-2}{4}}
\bigo{\ep^n \ex^{2t}} \di t \di
\theta^1 \ldots \di \theta^{n-1}\di z^1 \ldots \di z^{k} \\
\chi_1 S_{\gep}^R & = & R^{- \frac{n-2}{4}} \bigo{\ep^{-1}
\ex^{(n-1)t}} \\
\chi_2 S_{\gep}^R & = & R^{\frac{n-2}{4}} \bigo{\ep^{-1}
\ex^{-(n-1)t}} \\
\end{eqnarray*}
Moreover, by imposing that $\int_M F_\ep^R (v_\ep^R) \dvol{\gep}^R =
0 $ we have that
\begin{eqnarray*}
S^R & = & \frac{1}{\int_M 1 + \left( 1 + \frac{1}{c_n (n-1)} \right)
 v_\ep^R  + f(v_\ep^R) \dvol{\gep}^R }  \int_M S_{\gep}^R (1 +
v_\ep^R) \dvol{\gep}^R \\
& = & \frac{1}{\vol_{\gep}^R (M) }  \int_M S_{\gep}^R \dvol{\gep}^R
+ \littleo{\ep^{n-2}} \\
& = & R^{-\frac{n-2}{4}} \bigo{\ep^{n-2}}
\end{eqnarray*}
Therefore, for suitable constants $C_1, C_2$
\begin{eqnarray*}
\int_M \left(S^R - S^R_{\gep} \right) \chi_1 \dvol{\gep}^R & = & C_1
R^{\frac{m}{2} - \frac{n-2}{4}} \bigo{\ep^{n-2}}
\end{eqnarray*}
and
\begin{eqnarray*}
\int_M \left(S^R - S^R_{\gep} \right) \chi_2 \dvol{\gep}^R & = & C_2
R^{\frac{n-2}{4}} \bigo{\ep^{n-2}}
\end{eqnarray*}
Notice that $C_1$ and $C_2$ have the same sign, it is to say the
sign of both the terms $\int_M \left(S - S_{\gep} \right) \chi_1
\dvol{\gep}$ and $\int_M (S - S_{\gep} ) \chi_2 \dvol{\gep}$. Hence,
for large enough $R>0$ the sign of the approximate eigenvalue is
determined.\\

Obviously, if we make the same computation for the metric $Q g_2$
instead of $g_2$, we have that, for large enough $Q > 0$, the sign
of the approximate eigenvalue is the opposite of the one we found by
scaling the metric $g_1$.\\

Hence, if we look at the approximate eigenvalue as a continuous
function $\lambda_{F_\ep(v_\ep)} = \lambda_{F_\ep(v_\ep)} (P,Q)$
depending on the positive factors $P$ and $Q$, we have seen that
there exist $\bar{P}$ and $\bar{Q}$ such that
$\lambda_{F_\ep(v_\ep)}(\bar{P},1)$ and
$\lambda_{F_\ep(v_\ep)}(1,\bar{Q})$ have opposite sign.\\

At this point we can deduce that there exist positive real numbers
$P_0$ and $Q_0$ such that $\lambda_{F_\ep(v_\ep)} (P_0,Q_0) = 0$ and
this is what we wished.\\

\section{The non Ricci-flat case}

As we claimed in section 1, when both the initial metrics are scalar
flat but non Ricci flat it is possible to construct a zero scalar
curvature metric on the generalized connected sum. The idea consists
in doing a slight modification of the approximate solution metric
$g_\ep$ away from the gluing locus. By means of this modification it
is possible to obtain the vanishing of the term $\int_M F_\ep (v)
\dvol{\gep}$ without using the nonzero constant scalar curvature $S
= S(\ep,v)$ and it is also possible to show that up to carefully
choose the size of the adjustment, the approximate eigenvalue
$\lambda_{F_{\ep}(v_\ep)}$ is zero.\\

Let us describe the construction. Instead of the metric $g_\ep$ let
us consider the new approximate solution metric $\overline{\gep}
(r,s) = \gep + r h_1 + s h_2$, where $h_1$ and $h_2$ are positive
definite symmetric tensors supported respectively on the manifolds
$M_1$ and $M_2$ away from the polyneck, and $r$ and $s$ are real
numbers. Hence the equation we are induced to solve is the following
\begin{eqnarray}
\label{overline yamabe} \Delta_{\overline{\gep}} v & = &
\overline{F_\ep} (v,r,s) - \lambda_{\overline{F_\ep} (v,r,s)}
\beta_\ep
\end{eqnarray}
where $\overline{F_\ep} (v) = -c_n S_{\overline{\gep}} (1 + v)$.
Notice that by defining $\overline{F_\ep}$ that way we have
automatically imposed
that the final scalar curvature which we are going to achieve is zero.\\

As in the previous case, we will obtain the solution of the problem
by means of a fixed point argument lying on a sequence method.
Concerning the linear analysis, notice that the construction above
allows us to use all the results we have already obtained, hence let
us focus on the nonlinear analysis.\\

Since the condition $\int_M \overline{F_\ep} (v,r,s) \dvol{\gep} =
0$ has to be verified, we have that $r,s = \bigo{\ep^{n-2}}$ and so
it is easy to obtain a result analogous to the lemma \ref{Ball in
Ball} for the map $\overline{P_\ep}$ (with obvious notation), with
the same estimate for the radius $r_\ep$.\\

It is also immediate to prove that for sufficiently small $\ep>0$
the map $\overline{P_\ep}$ is a contraction and more precisely
\begin{eqnarray}
\norg{\overline{P_\ep}(u) - \overline{P_\ep}(v)}{M} & \leq & C \ep
\norg{u - v}{M}
\end{eqnarray}
for a suitable constant $C>0$. In particular the sequence defined by
$v_j := \overline{P_\ep^j (0)}$ converges with respect to the norm
$\norg{\cdot}{M}$ to a function $v_\ep$.\\

It remains to check that the operator $\overline{P_\ep}$ is actually
well defined. In particular we require that $\int_M \overline{F_\ep}
(v,r,s) \dvol{\gep} = 0$.\\

Before starting the calculation let us make some remarks concerning
the scalar curvature of the metric $\overline{\gep}$, in order to
get more information about $\overline{F_\ep}$. Since the supports of
$S_{\gep}, h_1, h_2$ are disjoint, we can write
\begin{eqnarray*}
S_{\overline{g_\ep}} & = & S_{\gep} + S_{g_1 + r h_1} + S_{g_2 + s
h_2} \\
& = & S_{\gep} \\
& + & S_{g_1} + r K_1 + \bigo{r^2} \\
& + & S_{g_2} + s K_2 + \bigo{s^2} \\
& = & S_{\gep} + r K_1 + \bigo{r^2} + s K_2 + \bigo{s^2}
\end{eqnarray*}
where, following \cite{Bes}
\begin{eqnarray*}
K_i & = & \Delta_{g_i} ({\rm tr}_{g_i} h_i) + \delta_{g_i}
(\delta_{g_i} h_i) + g_i ({\rm Ric}_{g_i}, h_i)
\end{eqnarray*}
for $i=1,2$. In the notation above $\delta_{g_i}$ indicate the
divergence of a symmetric tensor with respect to the metric $g_i$,
and ${\rm Ric}_{g_i}$ is the Ricci tensor of the metric $g_i$.\\

When we integrate we obtain
\begin{eqnarray*}
\int_{M} K_i \dvol{\gep} & = & \int_{M_i} K_i \dvol{g_i} \\
& = & \int_{M_i} g_i ({\rm Ric}_{g_i}, h_i) \dvol{g_i}
\end{eqnarray*}
because of the divergence theorem. Notice that in the Ricci flat
case the integral above is zero and there is no chance to correct
the term $\int_M S_{\gep}(1 + v) \dvol{\gep}$ in order to get the
condition
$\int_M \overline{F_\ep(v,r,s)} \dvol{\gep} = 0 $.\\

Let us define the map $G_{\ep,j} (r,s)$ as follows:
\begin{eqnarray*}
G_{\ep,j} & := & \int_M S_{\gep}(1 + v_j) \dvol{\gep}\\
& + & \int_{M_1} S_{g_1 + r h_1}(1 + v_j) \dvol{g_1}\\
& + & \int_{M_2} S_{g_2 + s h_2}(1 + v_j) \dvol{g_2}\\
& = & \int_{M} S_{\gep} \dvol{\gep} + r \int_{M_1} K_1 \dvol{g_1} +
s \int_{M_2} K_2 \dvol{g_2}\\
& + & E_1^{(j)}(r,s) + E_2^{(j)}(r,s)
\end{eqnarray*}
where
\begin{eqnarray*}
E_1^{(j)} & = & \int_{M} S_{\gep} v_j \dvol{\gep} + r \int_{M_1} K_1
v_j \dvol{g_1} + s \int_{M_2} K_2 v_j \dvol{g_2}\\
E_2^{(j)} & = & \int_{M_1} S_{g_1 + r h_1}(1 + v_j) \dvol{g_1} -
r \int_{M_1} K_1 (1 + v_j) \dvol{g_1}\\
& + & \int_{M_2} S_{g_2 + s h_2}(1 + v_j) \dvol{g_2} - s \int_{M_2}
K_2 (1 + v_j) \dvol{g_2}
\end{eqnarray*}
at this point our purpose is to describe the set where $G_{\ep,j}
(r,s)$ is zero.\\

Towards this aim let us consider the map $H_\ep (r,s) :=
G_{\ep,j}(r,s) - E_1^{(j)}(r,s) - E_2^{(j)}(r,s)$. In order to
simplify the calculus we can suppose that the symmetric tensors
$h_1$ and $h_2$ are so chosen that $\int_{M_1} K_1 \dvol{g_1} =
\int_{M_2} K_2 \dvol{g_2} = 1$. We can also assume $\int_{M}
S_{\gep} \dvol{\gep} < 0$ (if it is not the case, we can conclude by
means of obvious modifications) and since $\int_{M} S_{\gep}
\dvol{\gep} = \bigo{\ep^{n-2}}$ we can set, up to normalize,
$\int_{M} S_{\gep} \dvol{\gep} = - \ep^{n-2}$. The expression for
$H_\ep$ becomes then the following
\begin{eqnarray*}
H_\ep (r,s) & = & -\ep^{n-2} + r + s
\end{eqnarray*}
The set where $H_\ep$ vanishes is given by $\{ (r,s) \in
\mathbb{R}^2 \quad | \quad r + s = \ep^{n-2} \}$. We will show that
the set  where $G_{\ep,j}$ vanishes is uniformly close to the set
$\{ H_\ep = 0 \}$ with respect to
$j$.\\

As we have already seen, there exists positive constant $A>0$ such
that for every $j \in \mathbb{N}$
$$
\begin{array}{cccc}
  |v_j| & \leq & A \ep^{(n-2) + 1} & {\rm when} \quad n=3\\
  |v_j| & \leq & A \ep^{(n-2) + (1 + \gamma)} & {\rm when} \quad n \geq 4
\end{array}
$$
Moreover, since $r,s = \bigo{\ep^{n-2}}$ it is easy to see that
there exist positive constants $B_1, B_2$ such that for every $j \in
\mathbb{N}$
$$
\begin{array}{cccc}
|E_1^{(j)}| & \leq & B_1 \ep^{(n-2)+1} & {\rm when} \quad n=3\\
|E_1^{(j)}| & \leq & B_1 \ep^{(n-2) + (1 + \gamma)} & {\rm when} \quad n \geq 4\\
|E_2^{(j)}| & \leq & B_2 \ep^{2(n-2)} & \\
\end{array}
$$
In particular, for an arbitrarily small fixed constant $c>0$ and
sufficiently small $\ep>0$ we have that
\begin{eqnarray*}
|E_1^{(j)}| & \leq & (c/2) \ep^{n-2}\\
|E_2^{(j)}| & \leq & (c/2) \ep^{n-2}
\end{eqnarray*}
At this point it is immediate to see that for every $j \in
\mathbb{N}$
\begin{eqnarray*}
\{ G_{\ep,j} (r,s) = 0 \} & = & \{ (r,s) \in \mathbb{R}^2 \quad |
\quad r + s = \ep^{n-2} - E_1^{(j)}(r,s) - E_2^{(j)}(r,s) \}\\
& \subseteq & \{ (r,s) \in \mathbb{R}^2 \quad | \quad (1-c)
\ep^{n-2} \leq r + s \leq (1+c) \ep^{n-2} \}\\
& =: & Z_\ep
\end{eqnarray*}

If we set now $r_0^{(j)} := \ep^{n-2}/2$ for every $j \in
\mathbb{N}$, it must exists a real number $s_0^{(j)}$ such that
$(\ep^{n-2}, s_0^{(j)}) \in Z_\ep$ and $G_{\ep,j} (\ep^{n-2}/2,
s_0^{(j)}) = 0$.\\

Obviously $G_{\ep,j}$ are smooth functions with respect to the
variables $r$ and $s$ and in particular it is quite easy to get the
following uniform estimate for the first and the second partial
derivatives at the origin.
\begin{eqnarray*}
\left| \frac{\partial G_{\ep,j}}{\partial r} (0,0)  \right| & \leq &
\int_{M_1} \left| \frac{\partial S (g_1 + rh_1)}{\partial r}(0,0) \,
\cdot \, (1 + v_j) \right| \dvol{g_1} \\
& \leq & \nor{1 + v_j}{\mathcal{C}^0 (M)} \int_{M_1} \left|
\frac{\partial S (g_1 + rh_1)}{\partial r}(0,0) \right| \dvol{g_1} \\
& \leq & 2 \int_{M_1} \left| \frac{\partial S (g_1 + rh_1)}{\partial
r}(0,0) \right| \dvol{g_1} \\
\end{eqnarray*}
\begin{eqnarray*}
\left| \frac{\partial^2 G_{\ep,j}}{\partial r^2} (0,0)  \right| &
\leq & \int_{M_1} \left| \frac{\partial^2 S (g_1 + rh_1)}{\partial
r^2}(0,0) \,
\cdot \, (1 + v_j) \right| \dvol{g_1} \\
& \leq & 2 \int_{M_1} \left| \frac{\partial^2 S (g_1 +
rh_1)}{\partial r^2}(0,0) \right| \dvol{g_1} \\
\end{eqnarray*}
for every $j \in \mathbb{N}$. Notice that the same estimates hold
for the partial derivatives in $s$ and that the terms $
\frac{\partial^2 G_{\ep,j}}{\partial s \partial r} (0,0)$ and $
\frac{\partial^2 G_{\ep,j}}{\partial r \partial s} (0,0)$
are zero.\\

Concerning the modulus of the first partial derivative at the
origin, we are also able to provide a uniform lower bound, in fact
\begin{eqnarray*}
\left| \frac{\partial G_{\ep,j}}{\partial r} (0,0)  \right| & = &
\left| \int_{M_1} K_1 \, \cdot \, (1 + v_j)  \dvol{g_1} \right| \\
& \geq & \left| \int_{M_1} K_1 \dvol{g_1}  \right| - \int_{M_1}
|K_1| |v_j| \dvol{g_1} \\
& \geq & \left| \int_{M_1} K_1 \dvol{g_1}  \right| -
\nor{v_j}{\mathcal{C}^0 (M)}
\int_{M_1} |K_1| \dvol{g_1} \\
& \geq & \frac{1}{2} \left| \int_{M_1} K_1 \dvol{g_1}  \right|  \\
& > & 0 \quad {\rm for \,\, sufficiently \,\, small } \,\,\,\, \ep>0
\end{eqnarray*}
for every $j \in \mathbb{N}$. Of course, the same is true for
$\frac{\partial G_{\ep,j}}{\partial s} (0,0)$.

Now, arguing by contradiction and using these estimate it is
possible to deduce that there exists a positive constant $C>0$ and a
positive real number $R>0$ such that both the first partial
derivatives $\left|  \frac{\partial G_{\ep,j}}{\partial r}  \right|$
and $\left|  \frac{\partial G_{\ep,j}}{\partial r}  \right|$ are
greater than $C$ in $B_R \left( (0,0) \right)$, for every $j \in
\mathbb{N}$.\\

Up to choose $\ep$ sufficiently small, we have that the set $Z_\ep
\cap \{r,s \geq 0 \}$ lies in the ball of radius $R$ centered at the
origin, hence it is possible to apply the implicit function theorem
to the functions $G_{\ep,j}$ around the points $(\ep^{n-2}/2,
s_0^{(j)})$, so that we obtain, for every $j$, an open neighborhood
$U^{(j)}$ of $r_0^{(j)}$, an open neighborhood $V^{(j)}$ of
$s_0^{(j)}$ and a smooth function $f_j : U^{(j)} \longrightarrow
V^{(j)}$ such that $G_{\ep,j} (r, f_j(r)) = 0$ for every $r \in
U^{(j)}$.\\

Since it is possible to extend each implicit function $ f_j$ to the
interval $(0, (1-c)\ep^{n-2})$, we can suppose that there exists an
open neighborhood $U$ of $\ep^{n-2}/2$ and an open neighborhood $V$
of every $s_0^{(j)}$ such that it is possible to choose $U^{(j)} =
U$ and $V^{(j)} = V$ for every $j \in \mathbb{N}$.\\

Let us fucus now on the family of functions $\{f_j\}_{j \in
\mathbb{N}}$. Since each $f_j$ is a uniformly continuous function,
we can extend them to the compact set $\overline{U}$, so that we
have to handle now a family of functions $f_j : \overline{U}
\longrightarrow \overline{V}$ defined on a compact set and all
bounded by the same constant $(1+c) \ep^{n-2}$.\\

At this point our aim is to show that the $f_j$'s admit the same
Lipschitz's constant. First remember that
\begin{eqnarray*}
f_j(r) & = & \int_M S_{g_\ep} \dvol{g_\ep} - r + E_1^{(j)}(r) +
E_2^{(j)} (r)
\end{eqnarray*}
and consequently, for $r, r' \in U$
\begin{eqnarray*}
\left|f_j (r) - f_j (r')\right| & \leq & \left| r- r' \right| \\
& + & \left| E_1^{(j)}(r) - E_1^{(j)}(r') \right| \\
& + & \left| E_2^{(j)}(r) - E_2^{(j)}(r') \right|\\
& \leq & |r-r'|\\
& + & \int_{M_1} |K_1| \dvol{g_1} \nor{v_j}{\mathcal{C}^0(M)}
|r-r'|\\
& + & \int_{M_2} |K_2| \dvol{g_2} \nor{v_j}{\mathcal{C}^0(M)}
\left|f_j
(r) - f_j (r')\right|\\
& + & C_1 \ep^{n-2} |r-r'| + C_2 \ep^{n-2} \left|f_j (r) - f_j
(r')\right|
\end{eqnarray*}
for suitable $C_1,C_2>0$. Now, for $n=3$, we have that
$\nor{v_j}{\mathcal{C}^0(M)} \leq \tilde{D} \ep^{1-\gamma}$,
whereas, for $n \geq 4$, $\nor{v_j}{\mathcal{C}^0(M)} \leq \tilde{D}
\ep$, so it is easy to deduce that
\begin{eqnarray*}
\left| f_j(r) - f_j(r') \right| &\leq & \frac{1 + D_1 \ep^{1-\gamma}
+ C_1 \ep^{n-2}}{1 - D_2\ep^{1-\gamma} - C_2 \ep^{n-2}} |r - r'|
\end{eqnarray*}
when $n=3$, and
\begin{eqnarray*}
\left| f_j(r) - f_j(r') \right| &\leq & \frac{1 + D_1 \ep + C_1
\ep^{n-2}}{1 - D_2\ep - C_2 \ep^{n-2}} |r - r'|
\end{eqnarray*}
when $n \geq 4$, where $D_i = \tilde{D}\int_{M_i} |K_i| \dvol{g_i}$,
$i = 1,2$.\\

Thanks to the Ascoli - Arzel\'a theorem and up to consider
subsequence, the $f_j$'s converges with respect to the norm
$\nor{\cdot}{\mathcal{C}^0(U)}$ to a continuous function $f$.\\

Notice that in order to get a solution for the equation
\ref{overline yamabe} we could have picked a point $\overline{r} \in
\overline{U}$ and deduced immediately the convergence, up to
subsequences, of the $ f_j (\overline{r})$'s to a point
$\overline{s}$, but in the following we want to use the parameter
$r$ to kill the term $\lambda_{\overline{F_\ep}(v_\ep, r, f(r) )}$
and the
continuity of $f$ will be crucial.\\

We are now ready to discuss the sign of the term
$\lambda_{\overline{F_\ep}(v_\ep,r,f(r))}$ which appears in the
identity
\begin{eqnarray*}
\Delta_{\overline{\gep}} v_\ep & = & \overline{F_\ep} (v_\ep, r,
f(r)) - \lambda_{\overline{F_\ep} (v_\ep, r, f(r))} \beta_\ep
\end{eqnarray*}
As we have already seen, the sign of the approximate eigenvalue for
small $\ep>0$ is determined by the sign of its part of order
$\ep^{n-2}$, it is to say that if this one is strictly positive or
negative, so does $\lambda_{\overline{F_\ep}(v_\ep,r,f(r))}$. Hence
our task is to show that the sign of
$$\int_M S_{\gep}\chi_1 \dvol{\gep} - \int_M S_{\gep}\chi_2 \dvol{\gep}
+ r \int_{M_1} K_1 \dvol{g_1} - f(r) \int_{M_2} K_2 \dvol{g_2}$$
becomes positive or negative if we move $r\in \overline{U}$. Then,
thanks to the continuity of $f$ we deduce the existence of a real
number $\overline{r}$ in correspondence of which the approximate
eigenvalue vanishes.\\

Since there exists a constant $C>0$ such that
\begin{eqnarray*}
\left| \int_M S_{\gep}\chi_1 \dvol{\gep} - \int_M S_{\gep}\chi_2
\dvol{\gep} \right| & \leq & C\ep^{n-2}
\end{eqnarray*}
and since it is always possible to choose $r$ either in a region of
$U$ where $f(r)>r$ or in a region where $f(r)<r$, it is clear that
our goal is achieved if we impose that
$$\left|\int_{M_1} K_1 \dvol{g_1}\right| \,\, = \,\, \left| \int_{M_2} K_2 \dvol{g_2} \right|
\,\, > \,\, \frac{C \ep^{n-2}}{|f(r) - r|} $$ and we are always
allowed to do that.

\end{document}